\documentclass[11pt,jsco]{academic}
\usepackage{natbib}
\usepackage{amsmath}
\usepackage{amsfonts}
\usepackage{amssymb}
\usepackage{graphicx}%
\newenvironment{acknowledgement}{\noindent\textbf{Acknowledgements.}}{}

\newtheorem{theor}{Theorem}

\newtheorem{lem}{Lemma}
\newtheorem{obs}{Observation}

\def\Res{{\rm{Res}}}

\begin{document}


\title{An Elementary Proof of Sylvester's Double Sums for Subresultants}
\author{Carlos D' Andrea \thanks{Department d'\`Algebra i Geometria, Facultat de Matem\`atiques, Universitat
de Barcelona. Gran Via de les Corts Catalanes, 585; 08007 Spain.
\texttt{cdandrea@ub.edu}} , \;\; Hoon Hong
\thanks{Department of Mathematics, North Carolina State University,
Raleigh NC 27695, USA. Research supported by NSF CCR-0097976.
\texttt{hong@math.ncsu.edu}} , \;\; Teresa Krick
\thanks{Departamento de Matem\'atica, Facultad de Ciencias Exactas y
Naturales, Universidad de Buenos Aires, Ciudad Universitaria, 1428
Buenos Aires, Argentina. Research supported by grants CONICET PIP
2461/01 and UBACYT X-112. \texttt{krick@dm.uba.ar}} ,\;\; Agnes
Szanto \thanks{Department of Mathematics, North Carolina State
University, Raleigh NC 27695, USA. Research supported by NSF
grants CCR-0306406 and CCR-0347506. \texttt{aszanto@ncsu.edu} }}

\shortauthor{C. D'Andrea, H. Hong, T. Krick, A. Szanto}
\shorttitle{Sylvester's Double Sums}
\date{\today}
\maketitle


\begin{abstract}
In 1853 Sylvester stated and proved an elegant formula that expresses the
polynomial  subresultants in terms of the roots of the input
polynomials. Sylvester's formula was also recently proved  by
Lascoux and Pragacz by using multi-Schur functions and divided
differences. In this paper, we provide an elementary proof that uses only basic
properties of matrix multiplication and Vandermonde determinants.

\end{abstract}

\section{Introduction}

Subresultants play  a fundamental role in Computer Algebra and
Computational Algebraic Geometry (for instance, see
\cite{Collins:67,Brown_Traub:71,Collins:75,Renegar:92a,Gonzalez_Lombardi_Recio_Roy:89,
Gonzalez-Vega:96,Lombardi:2000,Hong:2000a,Apery_Jouanolou:05}). In
\cite{Sylvester:1853} Sylvester stated and proved an elegant formula that
expresses the polynomial subresultants of two polynomials in terms
of their roots, the so-called \emph{double-sum} formula. This
identity was proved also by Lascoux and Pragacz in
\cite{Lascoux_Pragacz:2001}, by using the theory of multi-Schur
functions and divided differences.

In this paper we provide a new and elementary  proof that uses
only the basic properties of matrix multiplication and Vandermonde
determinants.  As apparent in our proof,  Sylvester's double-sum
formula is only one simple step further a particular case,  the
so-called \emph{single-sum} formula. Such connection  between the
single and the double-sum formulae was originally thought to be
unlikely, as remarked in page 691 of \cite{Lascoux_Pragacz:2001}.
There have been various proofs for the single-sum formula
\cite{Apery_Jouanolou:05,Borchardt:60,Chardin:90,Hong:99a,Diaz-Toca:2004}.

The matrix multiplication technique, presented in this papers, has proven to be quite powerful
in that it is easily generalizable  to multivariate polynomials:
similar techniques were successfully applied to obtain expressions
for multivariate subresultants in roots in \cite{DKS05}, and the
generalization of Sylvester's single and double-sum formulae to
the multivariate case is the subject of ongoing research.

\section{Review of Sylvester's Double Sum for Subresultants}

Let $f=a_{m}x^{m}+\cdots+a_{0}$ and $g=b_{n}x^{n}+\cdots+b_{0},$ be
two polynomials with coefficients in a commutative ring. The $d$-th
\emph{subresultant polynomial}
$\operatorname*{Sres}\nolimits_{d}(f,g)$ is defined for $0\le d<
\min\{m,n\}$  or, if $m\neq n$ holds,  for   $d=\min\{m,n\}$, as the
following determinant:
\begin{equation}\label{defs}
\operatorname*{Sres}\nolimits_{d}(f,g)
:=\det%
\begin{array}{|cccccc|c}
\multicolumn{6}{c}{\scriptstyle{m+n-2d}}\\
\cline{1-6}
a_{m} & \cdots & & \cdots & a_{d+1-\left(n-d-1\right)}& x^{n-d-1}f(x)&\\
& \ddots & && \vdots  & \vdots &\scriptstyle{n-d}\\
&  &a_{m}&\cdots &a_{d+1}&f(x)& \\
\cline{1-6}
b_{n} &\cdots & & \cdots & b_{d+1-(m-d-1)}&x^{m-d-1}g(x)&\\
&\ddots &&&\vdots & \vdots &\scriptstyle{m-d}\\
&& b_{n} &\cdots & b_{d+1} & g(x)&\\
\cline{1-6}
\multicolumn{2}{c}{}
\end{array}
\end{equation}
where $a_{\ell}=b_{\ell}=0$ for $\ell<0$.
\par
By developing this determinant by the last column, it is clear that
$\operatorname*{Sres}\nolimits_{d}(f,g)$ is a polynomial combination
of $f$ and $g$. It is also a classic fact that
$\operatorname*{Sres}\nolimits_{d}(f,g)$ is a polynomial of degree
bounded by $d$, since it coincides with the determinant of the
matrix obtained by replacing the last column $C_{m+n-2d}$ by
\begin{center}$C'_{m+n-2d}:= C_{m+n-2d}-x^{d+1}C_{m+n-2d-1} -\cdots -
x^{m+n-d-1}C_1$.\end{center}

Now, let $A=(\ldots,\alpha,\ldots)$ and $B=(\ldots,\beta,\ldots)$
be finite lists (ordered sets) of distinct indeterminates. In
\cite{Sylvester:1853} Sylvester introduced for $0\le p\le |A|,
0\le q\le |B|$ the following {\em double-sum} expression in $A$
and $B$:

\[
\operatorname*{Sylv}\nolimits^{p,q}(A,B;x)  :=
\sum_{\substack{A^{\prime }\subset A,\,B^{\prime}\subset
B\\|A^{\prime}|=p,\,|B^{\prime}|=q}}R(x,A^{\prime
})\,R(x,B^{\prime})\,\frac{R(A^{\prime},B^{\prime})\,R(A\backslash
A^{\prime},B\backslash B^{\prime})}{R(A^{\prime},A\backslash
A^{\prime })\,R(B^{\prime},B\backslash B^{\prime})},\] where
\[
R(X,Y):=\prod_{x\in X,y\in Y}(x-y), \ \quad R(x,Y):=\prod_{y\in
Y}(x-y).
\]

In \cite{Sylvester:1853}  Sylvester gave the following  elegant
formula that expresses the subresultants in terms of the double-sum,
that is, in terms of the roots of $f$ and $g$.

\begin{theor}
[Sylvester's double-sum formula]\label{sylvester} Let  $f, g$ be the
monic polynomials
\[
f=\prod_{\alpha\in A}(x-\alpha),\;\;\;g=\prod_{\beta\in B}(x-\beta
)\ \ \in\mathbb{Z}[\alpha\in A,\,\beta\in B][x],
\]
 where $|A|=m$ and $|B|=n$. Let $p,q\ge 0$ be such that
$d:=p+q<\min\{m,n\}$ or $d=\min\{m,n\}$ if $m\neq n$ holds. Then
$$
\operatorname*{Sres}\nolimits_{d}(f,g)=\frac{(-1)^{p(m-d)}}
{\binom{d}{p}}\;\operatorname*{Sylv}\nolimits^{p,q}%
(A,B;x).
$$
\end{theor}

\noindent When $p=d$ and $q=0$, the above expression immediately simplifies to
the $\emph{single}$-$\emph{sum\ }$formula:%
\begin{equation}\label{single-sum formula}
\operatorname*{Sres}\nolimits_{d}(f,g)=\sum_{\substack{A^{\prime}\subset
A\\|A^{\prime}|=d}}R(x,A^{\prime})\;\frac{R(A\backslash A^{\prime}%
,B)}{R(A\backslash A^{\prime},A^{\prime})}.%
\end{equation}

Complete proofs of Sylvester's double-sum can be found in
\cite{Sylvester:1853,Lascoux_Pragacz:2001}, while the single-sum formula has
various proofs,
\cite{Apery_Jouanolou:05,Borchardt:60,Chardin:90,Hong:99a,Diaz-Toca:2004}.
Here we present in Section~\ref{4} an alternative elementary proof
for both results.

\section{Notations}\label{3}

We recall that  $0\le d<\min\{m,n\}$ or $d:=\min\{m,n\}$ if $m\ne n$
holds.
We let $M_{f}$ and $M_{g}$ denote the following matrices: {\small
\[\begin{array}{ccc}
M_{f}:=
\begin{array}{|ccccc|c}
\multicolumn{5}{c}{\scriptstyle{m+n-d}}\\
\cline{1-5}
a_{0} & \dots & a_{m} &  & &\\
& \ddots &  & \ddots & &\scriptstyle{n-d}\\
&  & a_{0} & \dots & a_{m} & \\
\cline{1-5}
\multicolumn{2}{c}{}
\end{array},
&
&
M_{g}:=
\begin{array}{|ccccc|c}
\multicolumn{5}{c}{\scriptstyle{m+n-d}}\\
\cline{1-5}
b_{0} & \dots & b_{n} &  & &\\
& \ddots &  & \ddots & &\scriptstyle{m-d}\\
&  & b_{0} & \dots & b_{n} & \\
\cline{1-5}
\multicolumn{2}{c}{}
\end{array}
\end{array}.
\]}
We now define {\small $$S_{d}:=%
\begin{array}{|c|c}
\multicolumn{1}{c}{ \scriptstyle{m+n-d}}&\\
\cline{1-1}
\ \ M_{t-x}\ \ &\scriptstyle{d}\\
\cline{1-1}
M_f&\scriptstyle{n-d}\\
\cline{1-1}
M_g&\scriptstyle{m-d}\\
\cline{1-1} \multicolumn{2}{c}{}
\end{array} \quad \mbox{where} \quad M_{t-x}:=
\begin{array}{|ccccccc|c}
\multicolumn{7}{c}{ \scriptstyle{m+n-d}}\\
\cline{1-7}
-x&1&0& \dots& &\dots&0 &\\
& \ddots &\ddots  & \ddots &&&\vdots &\scriptstyle{d}\\
 & & -x & 1&0&\dots & 0 & \\
\cline{1-7} \multicolumn{2}{c}{}
\end{array} .$$}
Finally, we define for a polynomial $p(t)$  and two lists,
$\Gamma:=(\gamma_1,\ldots,\gamma_u)$ of scalars  and $E:=(e_1,\dots,
e_v)$ of non-negative integers,  the (not-necessarily square) matrix
of size $v\times u$:
$$\langle p(t),\Gamma \rangle_E:=\begin{array}{|ccc|c}
\multicolumn{3}{c}{\scriptstyle{u}}&\\
\cline{1-3}
\gamma_1^{e_1}p(\gamma_1)&\dots & \gamma_u^{e_1}p(\gamma_u)&\\
\vdots& & \vdots& \scriptstyle v\\
\gamma_1^{e_v}p(\gamma_1)&\dots & \gamma_u^{e_v}p(\gamma_u)& \\
 \cline{1-3}
 \multicolumn{2}{c}{}
\end{array}.$$
For instance, under  this notation, if we take $E:=(0,\dots,u-1)$,
we have the following equality for the Vandermonde determinant
$\mathcal{V}(\Gamma)$ associated to $\Gamma$:
$$\mathcal{V}(\Gamma):=|(\gamma_j^{i-1})_{1\le i,j\le u}|=|\langle 1, \Gamma\rangle_E|.$$
 When
 $E$ is of the form $E=(0,\dots,v-1)$, we directly write $\langle
p(t),\Gamma \rangle_v$. \\We mention   the following useful
equalities that hold since $m+n-d\ge \max(m,n)$:
$$\begin{array}{lcl} M_f\cdot \langle 1,
\Gamma\rangle_{m+n-d}&=&\langle f(t),\Gamma\rangle_{n-d}\\
 M_g\cdot
\langle 1, \Gamma\rangle_{m+n-d}&=&\langle g(t),\
\Gamma\rangle_{m-d}\\ M_{t-x}\cdot \langle 1 ,
\Gamma\rangle_{m+n-d}&=&\langle t-x,\Gamma\rangle_{d}\end{array}.$$

\section{The Proof}\label{4}

The proof is divided into a series of lemmas which are interesting
on their own. For an easier understanding, we recommend not to pay
attention to signs in a first approach.

\begin{lem}\label{kd} Under the previous assumptions and notations,
we have
$$\operatorname*{Sres}\nolimits_{d}(f,g)=(-1)^{d+(n-d)(m-d)}|S_d|.$$
\end{lem}

\begin{proof}
We denote by $C_i$ the $i$-th column of the matrix $S_d$ and we
replace its first column $C_1$  by $C'_1:=C_1+xC_2+\ldots
+x^{m+n-d-1}C_{m+n-d}$. This operation does not change the
determinant of this matrix, and {\scriptsize $${\large C'_1}:=
\begin{array}{|c|c} \cline{1-1}
0&\\
\vdots &\scriptstyle{d}\\
0 & \\
\cline{1-1}
f(x)&\\
\vdots&\scriptstyle{n-d}\\
x^{n-d-1}f(x)&\\
\cline{1-1}
g(x)&\\
\vdots&\scriptstyle{m-d}\\
x^{m-d-1}g(x)&\\
\cline{1-1}
\end{array}.$$}
We now perform  a Laplace expansion of the determinant of the new
matrix over the first $d$ rows, and we observe that only one block
survives, which corresponds to columns $2$ to $d+1$ of $M_{t-x}$.
Moreover, this block is  lower triangular with diagonal entries $1$.
Thus {\small
\begin{eqnarray*}|S_d|&=&(-1)^d \det
\begin{array}{|cccccc|c}
\multicolumn{6}{c}{\scriptstyle{m+n-2d}}\\
\cline{1-6}
f(x) & a_{d+1} & \dots& a_m & & &\\
\vdots & \vdots & && \ddots  &  &\scriptstyle{n-d}\\
x^{n-d-1}f(x)& a_{d+1-(n-d-1)} &\dots& &\dots &a_m& \\
\cline{1-6}
g(x) & b_{d+1} & \dots& b_n & & &\\
\vdots & \vdots & && \ddots  &  &\scriptstyle{m-d}\\
x^{m-d-1}g(x)& b_{d+1-(m-d-1)} &\dots& &\dots &b_n& \\
\cline{1-6} \multicolumn{2}{c}{}
\end{array}\\
&=&(-1)^{d+(n-d)(m-d)}\operatorname*{Sres}\nolimits_{d}(f,g),
\end{eqnarray*}}
since the matrix in the right-hand side above is the matrix of
(\ref{defs}) viewed backward.
\end{proof}

For simplicity, from now on,  we assume  $f$ and $g$ to be the monic
polynomials $f=\prod_{\alpha\in A}(x-\alpha)$, $g=\prod_{\beta\in
B}(x-\beta )$ where $A$ and $B$ are lists with $|A|=m$ and $|B|=n$.
(As pointed out by a referee,  under this assumption   one has  in
the language of multi-Schur functions:
$|S_d|=S_{1^d;(m-d)^{n-d};0^{m-d}}(-x,-A,-B)$ (see
\cite{Lascoux_Pragacz:2001}).)

\medskip
The lemmas below generalize in an obvious manner to non-monic
polynomials. The first one corresponds to Th.~3 in \cite{Hong:99a}.
We prove it here with a different technique that follows from Lemma
\ref{kd}.

\begin{lem}
\label{MsVA_Det} \  {\it (Hong's subresultant in roots
\cite[Th.~3.1]{Hong:99a})}\\
Under the previous notations, we have
$$\operatorname*{Sres}\nolimits_{d}(f,g) \,
\mathcal{V}\left(  A\right) =\det
\begin{array}{|c|c}
\multicolumn{1}{c}{ \scriptstyle{m}}&\\[1mm]
\cline{1-1}
\langle x-t,A\rangle_d &\scriptstyle{d}\\[1mm]
\cline{1-1}
\langle g(t),A\rangle_{m-d}&\scriptstyle{m-d}\\[1mm]
\cline{1-1} \multicolumn{2}{c}{}
\end{array}\ .
$$
\end{lem}

\begin{proof}
We note that  $\left\vert S_{d}\right\vert \, \mathcal{V}(A)$ is the
determinant of the following  product of matrices:

$$
\begin{array}{c|c|}
\multicolumn{1}{c}{}&
\multicolumn{1}{c}{\scriptstyle{m+n-d}}\\
\cline{2-2}
\scriptstyle{d}& \ \ M_{t-x}\ \ \\
\cline{2-2}
\scriptstyle{n-d}& M_f\\
\cline{2-2} \scriptstyle{m-d}&
M_g\\
\cline{2-2} \multicolumn{2}{c}{}
\end{array}
\,
\begin{array}{|c|c|c}
\multicolumn{1}{c}{ \scriptstyle{m}}&\multicolumn{1}{c}{ \scriptstyle{n-d}}&\\
\cline{1-2}
& 0 &\scriptstyle{m}\\
\langle 1,A\rangle_{m+n-d}&\\
&I_{n-d} &\scriptstyle{n-d}\\
\cline{1-2} \multicolumn{3}{c}{}
\end{array}= \begin{array}{|c|c|c}
\multicolumn{1}{c}{ \scriptstyle{m}}&\multicolumn{1}{c}{ \scriptstyle{n-d}}&\\
\cline{1-2}
 \langle t-x,A\rangle_{d}& * &\scriptstyle{d}\\
\cline{1-2} \ 0\ &M'_f&\scriptstyle{n-d}\\
\cline{1-2}\langle g(t),A\rangle_{m-d}
& *  &\scriptstyle{m-d}\\
\cline{1-2} \multicolumn{3}{c}{}
\end{array},
$$
since $\langle f(t),
A\rangle_{n-d}=\left[\alpha_j^{i-1}f(\alpha_j)\right]=[0]$.\\
By permuting the rows of the second block with those of the third,
we obtain \begin{eqnarray*}
\operatorname*{Sres}\nolimits_{d}(f,g)\,\mathcal{V}(A)&=&(-1)^{d+(m-d)(n-d)}|S_{d}|\,\mathcal{V}(A)
\\&=& (-1)^d \det
\begin{array}{|c|}
\cline{1-1}
 \langle t-x,A\rangle_{d}\\
\cline{1-1}
\langle g(t),A\rangle_{m-d} \\
\cline{1-1}
\end{array}\,|M_{f}^{\prime}|\\
&=& \det
\begin{array}{|c|}
\cline{1-1}
 \langle x-t,A\rangle_{d}\\
\cline{1-1}
\langle g(t),A\rangle_{m-d} \\
\cline{1-1}
\end{array},
\end{eqnarray*}
 since   $M_{f}^{\prime}$ is a
lower triangular matrix with diagonal entries $a_{m}=1$.
\end{proof}

Let us remark here that the Poisson product formula $\Res(f,g)=
\prod_{\alpha\in A}g(\alpha)$ is a direct consequence of the
previous Lemma for the case $d=0$.

\medskip
 For $S\subseteq T$ finite lists, let $\operatorname*{sg}\left( S,T\right)
$\thinspace$:=\left( -1\right) ^{\sigma}$ where $\sigma$ is the
number of transpositions needed to take $T$ to $S  \cup (T\backslash
S).$ Here, ``$\cup$" stands for list concatenation and
``$\backslash$" means list subtraction.

\begin{lem}
\label{MsVAVB_Det}  Let $P$ and $Q$ be two disjoint sublists of
$E:=(0,\dots,d-1)$ that satisfy $P\cup Q=E$, and let $p:=|P|$,
$q:=|Q|$. Then {\small
\begin{equation}\label{lema3}
\operatorname*{Sres}\nolimits_{d}(f,g)\,
\mathcal{V}(A)\mathcal{V}(B)  =(-1)^{q+(m-d)n}%
\operatorname*{sg}\left(  P,E\right)  \det
\begin{array}{|c|c|l}
\multicolumn{1}{c}{\scriptstyle m}&\multicolumn{1}{c}{\scriptstyle n}&\\
\cline{1-2}
\langle x-t, A\rangle_P&0&\scriptstyle p\\
\cline{1-2}
0&\langle x-t, B\rangle_Q&\scriptstyle q\\
\cline{1-2}
\langle 1, A\rangle_{m+n-d}&\langle 1, B\rangle_{m+n-d}&\scriptstyle{m+n-d}\\
\cline{1-2} \multicolumn{1}{c}{}
\end{array}.
\end{equation}}
\end{lem}

\begin{proof}
Recalling that $\mathcal{V}(B)=|\langle 1,B\rangle_{n}|$,   we have by Lemma \ref{MsVA_Det}:%
\small{\begin{eqnarray*}\operatorname*{Sres}\nolimits_{d}(f,g)\,\mathcal{V}(A)\mathcal{V}(B)&=&
\det
\begin{array}{|c|c|l}
\multicolumn{1}{c}{\scriptstyle{m}}&\multicolumn{1}{c}{\scriptstyle{n}}&\\
\cline{1-2}
\langle x-t,A\rangle_{d}&0&\scriptstyle{d}\\
\cline{1-2}
\langle g(t),A\rangle_{m-d}&0&\scriptstyle{m-d}\\
\cline{1-2}
\langle 1,A\rangle_{n}&\langle 1,B\rangle_{n}&\scriptstyle{n}\\
\cline{1-2} \multicolumn{1}{c}{}
\end{array}
\\ &=&
(-1)^{(m-d)n} \det
\begin{array}{|c|c|l}
\multicolumn{1}{c}{\scriptstyle{m}}&\multicolumn{1}{c}{\scriptstyle{n}}&\\
\cline{1-2}
\langle x-t,A\rangle_{d}&0&\scriptstyle{d}\\
\cline{1-2}
\langle 1,A\rangle_{n}&\langle 1,B\rangle_{n}&\scriptstyle{n}\\
\cline{1-2}
\langle g(t),A\rangle_{m-d}&0&\scriptstyle{m-d}\\
\cline{1-2} \multicolumn{1}{c}{}
\end{array}
\end{eqnarray*}
$$ =(-1)^{(m-d)n}\det\left(
\begin{tabular}
[c]{cccc}
& $\scriptstyle d$ & $\scriptstyle n$ & $\scriptstyle{m-d}$\\\cline{2-4}%
$\scriptstyle d$ & \multicolumn{1}{|c}{$I_d$} &
\multicolumn{1}{|c}{$0$} &
\multicolumn{1}{|c|}{$0$}\\\cline{2-4}%
$\scriptstyle n$ & \multicolumn{1}{|c}{$0$} &
\multicolumn{1}{|c}{$I_n$} &
\multicolumn{1}{|c|}{$0$}\\\cline{2-4}%
$\scriptstyle{m-d}$ & \multicolumn{1}{|c}{$0$} & \multicolumn{2}{|c|}{$M_{g}$}\\\cline{2-4}%
\end{tabular}
\,\
\begin{tabular}
[c]{ccc}%
$\scriptstyle m$ & $\scriptstyle n$ & \\\cline{1-2}%
\multicolumn{1}{|c}{$\langle x-t,A\rangle_d $} &
\multicolumn{1}{|c}{$0$} &
\multicolumn{1}{|l}{$\scriptstyle d$}\\\cline{1-2}%
\multicolumn{1}{|c}{\vspace{-3mm}} & \multicolumn{1}{|c}{ }
& \multicolumn{1}{|l}{}\\
\multicolumn{1}{|c}{\vspace{-2mm}$\langle 1,A\rangle_{m+n-d} $} &
\multicolumn{1}{|c}{$\langle 1,B\rangle_{m+n-d} $} & \multicolumn{1}{|l}{$\scriptstyle{m+n-d}$}\\
\multicolumn{1}{|c}{} & \multicolumn{1}{|c}{ } &
\multicolumn{1}{|l}{}\\
\cline{1-2}%
\end{tabular}
\right) ,
$$}
since $M_g\cdot \langle 1,B\rangle_{m+n-d}=\langle
g(t),B\rangle_{m-d}=[0]$. Now, since the first matrix is lower
triangular with diagonal entries $1$, we
have%
\begin{equation}\label{single-sum}
\operatorname*{Sres}\nolimits_{d}(f,g)\,\mathcal{V}(A)  \,
\mathcal{V}(B) =(-1)^{(m-d)n}\,\det
\begin{array}{|c|c|l}
\multicolumn{1}{c}{\scriptstyle m}&\multicolumn{1}{c}{\scriptstyle n}&\\
\cline{1-2}
\langle x-t,A\rangle_d&0&\scriptstyle d\\
\cline{1-2}
\langle 1,A\rangle_{m+n-d}&\langle 1,B\rangle_{m+n-d}&\scriptstyle{m+n-d}\\
\cline{1-2} \multicolumn{1}{c}{}
\end{array}
.\end{equation} Finally, recalling that $\langle
x-t,A\rangle_d=\left( \alpha_j^{i-1}x-\alpha_j^{i}\right)_{1\le i\le
d, 1\le j\le m}$ and\\
$\langle 1,A\rangle_{m+n-d}= \left( \alpha_j^{i-1}\right)_{1\le i\le
m+n-d, 1\le j\le m}$, the obvious subtractions and permutations of
rows yield%
{\small $$
\operatorname*{Sres}\nolimits_{d}(f,g)\,\mathcal{V}(A)\mathcal{V}(B)  =(-1)^{(m-d)n}%
\operatorname*{sg}\left(  P,E\right) \det
\begin{array}{|c|c|l}
\multicolumn{1}{c}{\scriptstyle m}&\multicolumn{1}{c}{\scriptstyle n}&\\
\cline{1-2}
\langle x-t,A\rangle_P&0&\scriptstyle p\\
\cline{1-2}
0&-\langle x-t,B\rangle_Q&\scriptstyle q\\
\cline{1-2}
\langle 1,A\rangle_{m+n-d}&\langle 1,B\rangle_{m+n-d}&\scriptstyle{m+n-d}\\
\cline{1-2} \multicolumn{1}{c}{}
\end{array}.
$$}
The lemma follows by moving $(-1)^q$ out of the determinant.
\end{proof}

We will also need in the proof the  following observation:
\begin{obs} Let $\Gamma:=(\gamma_1,\dots,\gamma_d)$. Then
\begin{equation}\label{tech}
|\langle x-t, \Gamma\rangle_d|=  R(x,\Gamma)|\langle
1,\Gamma\rangle_d|.\end{equation}
\end{obs}

\begin{proof} The claim follows from {\scriptsize{
$$\left(
\begin{array}{ccc}x-\gamma_1&\dots &
x-\gamma_d\\
\vdots &  & \vdots\\
\gamma_1^{d-1}x-\gamma_1^{d} & \dots &
\gamma_d^{d-1}x-\gamma_d^{d}\end{array}\right)=\left(\begin{array}{ccc}1&\dots&1\\
\vdots& & \vdots\\
\gamma_1^{d-1}&\dots&\gamma_d^{d-1}\end{array}\right)\,\left(\begin{array}{ccc}x-\gamma_1& & \\
 & \ddots &  \\
& &x-\gamma_d\end{array}\right).$$}}
\end{proof}

\subsection{Proof of Theorem \ref{sylvester}}
For any $P$ and $Q$  disjoint sublists of $E:=(0,\dots,d-1)$ that
satisfy $P\cup Q=E$,   with $| P|=p $ and $| Q|=q$,   a Laplace
expansion over the first $d$ rows in Identity (\ref{lema3}) gives
that  $\operatorname*{Sres}\nolimits_{d}(f,g)\,\mathcal{V}\left(
A\right) \mathcal{V}\left( B\right)$ equals {\small
$$
(-1)^{\sigma}\operatorname*{sg}( P,E) \!\!\!
\sum_{\substack{A^{\prime}\subset A,\,B^{\prime}\subset
B\\\left\vert A^{\prime}\right\vert =p,\left\vert
B^{\prime}\right\vert
=q}}\!\operatorname*{sg}(A^{\prime},A)\,\operatorname*{sg}(B^{\prime
},B)\cdot|\langle x-t,A'\rangle_P|\cdot| \langle
x-t,B'\rangle_Q|\cdot\mathcal{V}\left(A\backslash A^{\prime}\cup
B\backslash B^{\prime}\right)
$$
} where $\sigma:=q+(m-d)n+(m-p)q \equiv (m-d)(n-q) \pmod 2$. Adding
over all such choices of $P\subset E$ with $|P|=p$, we deduce that
$\operatorname*{Sres}\nolimits_{d}(f,g)\,\mathcal{V}\left( A\right)
\mathcal{V}\left( B\right)$ equals {\small
  $$
  \frac{1}{\binom{d}{p}} \sum_{P}\left(
-1\right) ^{\sigma}\operatorname*{sg}\left( P,E\right)
\sum_{A^\prime,B^\prime}\,\operatorname*{sg}(A^{\prime},A)
\operatorname*{sg}(B^{\prime },B)\cdot|\langle
x-t,A'\rangle_P|\cdot| \langle
x-t,B'\rangle_Q|\cdot\mathcal{V}\left(A\backslash A^{\prime}\cup
B\backslash B^{\prime}\right)$$ $$= \
\frac{(-1)^{\sigma}}{\binom{d}{p}}\!\!
\sum_{A^\prime,B^\prime}\!\operatorname*{sg}(A^{\prime},A)\!\operatorname*{sg}(B^{\prime},B)\mathcal{V}\left(
A\backslash A^{\prime} \cup\,B\backslash B^{\prime}\right)
\!\!\left( \sum_{P}\!\operatorname*{sg}\left(  P,E\right)  |\langle
x-t,A'\rangle_P|\cdot| \langle x-t,B'\rangle_Q| \right). $$} We
observe now that, by another Laplace expansion and Identity
(\ref{tech}), {\small
$$\sum_{P}\!\operatorname*{sg}\left(  P,E\right)  |\langle
x-t,A'\rangle_P|\cdot| \langle x-t,B'\rangle_Q|=|\langle x-t,A'\cup
B'\rangle_d|=R(x,A')R(x,B')|\langle 1,A'\cup B'\rangle_d|.$$}
Recalling that  $|\langle 1,A'\cup
B'\rangle_d|=\mathcal{V}\left(A'\cup B'\right)$, this gives

{\small
\begin{eqnarray*}
\operatorname*{Sres}\nolimits_{d}(f,g)\,&=&
\frac{(-1)^{\sigma}}{\binom{d}{p}}
\sum_{A^\prime,B^\prime}\!R(x,A')R(x,B')\frac{\operatorname*{sg}(A^{\prime},A)\,\!\operatorname*{sg}(B^{\prime},B)\mathcal{V}\left(
A\backslash A^{\prime} \cup\,B\backslash
B^{\prime}\right)\mathcal{V}\left(A'\cup
B'\right)}{\mathcal{V}(A)\mathcal{V}(B)}\\
&=&\frac{(-1)^{\sigma}(-1)^{\tau}}{\binom{d}{p}}
\sum_{A^\prime,B^\prime} R(x,A')R(x,B')\frac{R(A',B')R(A\backslash
A',B\backslash B')}{R(A',A\backslash A')R(B',B\setminus B')},
\end{eqnarray*}}
where $\tau= (m-p)(n-q)+pq-(m-p)p-(n-q)q=(m-d)(n-d)$ since for any
finite lists $X,Y$, one has $\mathcal{V}(X\cup
Y)=\mathcal{V}(X)\mathcal{V}(Y)R(Y,X)=(-1)^{|X|\cdot|Y|}\mathcal{V}(X)\mathcal{V}(Y)R(X,Y)$.
\\
The claim follows now  from the fact that
$(m-d)(n-q)+(m-d)(n-d)\equiv  (m-d)p \pmod 2$.{\hfill\mbox{$\Box$}

\bigskip
As a final remark, we mention that  if in the previous  proof we
start with  a Laplace expansion over the first $d$ rows in Identity
(\ref{single-sum}) instead of Identity  (\ref{lema3}), we obtain in
the same manner  Sylvester's single sum formulation (\ref{single-sum
formula}).

\bigskip

\begin{acknowledgement}
We are grateful to the anonymous referees for their careful reading of
 our preliminary manuscript and their very precise indications to  improve our presentation.
\end{acknowledgement}


\end{document}